\documentclass [11pt]{article}
\usepackage{amsfonts}
\usepackage{amsmath}
\newtheorem{theo}{Theorem}[section]
\newtheorem{prop}[theo]{Proposition}
\newtheorem{lemma}[theo]{Lemma}

\title{On the largest eigenvalue of a   random subgraph of the hypercube
}
\author{Alexander Soshnikov\thanks{
Department of Mathematics,
University of California at Davis, 
One Shields Ave., Davis, CA 95616, USA.
Email address: soshniko@math.ucdavis.edu.
Research was supported in part by the NSF grant DMS-0103948.
} 
\and Benny Sudakov \thanks{Department of Mathematics,
Princeton University, Princeton, NJ 08544, USA
and Institute for Advanced Study, Princeton, NJ 08540,
USA. Email address: bsudakov@math.princeton.edu.
Research was supported in part by NSF grants DMS-0106589, CCR-9987845
and by the State of New Jersey.
}  }
\date{}
\begin{document}
\maketitle
\begin{abstract}
Let $G$ be a random subgraph of the $n$-cube where each edge appears
randomly and independently with  probability $p$. We prove that 
the largest eigenvalue of the adjacency matrix of $G$ is almost surely
$$\lambda_1(G)=(1+o(1)) \max \big( \Delta^{1/2}(G), n\*p \big),$$ 
where $\Delta(G)$ is the maximum degree of $G$ and the $o(1)$ term tends to 
zero as $\max( \Delta^{1/2}(G), n\*p )$ tends to infinity.
\end{abstract}

\section{Introduction and formulation of results}
Let $ Q^n$ be a graph whose vertices are all the 
vectors $\big\{ x=(x_1,\ldots, x_n)~|~x_i \in \{0,1\}\big\} $ and two vectors 
$x$ and $y$ are 
adjacent if they differ in exactly one coordinate, i.e., 
$\sum_i |x_i-y_i|=1$. We call $Q^n$ the {\em $n$-dimensional cube} or simply
the {\em $n$-cube}. Clearly $Q^n$ is an $n$-regular, bipartite graph of 
order $2^n$.
In this paper we study random subgraphs of the $n$-cube.
A random subgraph $G(Q^n,p)$ is a discrete probability space composed of 
all subgraphs of $n$-cube, where each edge of $Q^n$ appears randomly and 
independently with probability $p=p(n)$.
Sometimes with some abuse of notation we will refer to a random subgraph
$G(Q^n,p)$ as a graph on $2^n$ vertices generated according to the
distribution described above.
Usually asymptotic properties of random graphs are of interest. We say
that a graph property ${\cal A}$ holds {\em almost surely}, or a.s. for
brevity, in $G(Q^n,p)$ if the probability that $G(Q^n,p)$ has property 
${\cal A}$
tends to one as $n$ tends to infinity. Necessary
background information on random graphs in general and random subgraphs of
$n$-cube can be found in [4].

Random subgraphs of the hypercube were  studied by Burtin [5], Erd\"os and 
Spencer [8], Ajtai, Koml\'os and Szemer\'edi  [2] and 
Bollob\'as [4], among 
others. In particular it was shown that a giant component emerges shortly after
$p=1/n$ ([2]) and the graph becomes connected shortly after 
$p=1/2$ ([5], [8], [4]).
Recently the model has become of interest in mathematical 
biology ([7], [13], [14]). In particular it appears ( see [13],[14]) that 
random graphs play an important role in a general model of a population 
evolving over a network of selectively neutral genotypes. It has been shown
that the population's limit distribution on the neutral network is solely 
determined by the network topology and given by the principal eigenvector of 
the network's adjacency matrix. Moreover, the average number of neutral mutant 
neighbors per individual is given by the spectral radius. 

The subject of this
paper is the asymptotic behavior of 
the largest eigenvalue  of the random graph $G(Q^n, p)$.
The adjacency matrix of $G$  is an $2^n \times 2^n$ matrix $A$ whose 
entries are  either one or zero depending on whether the edge $(x,y)$ 
is present in $G$ or not. $A$ is a random real symmetric matrix with the 
eigenvalues denoted by $ \lambda_1 \geq \lambda_2 \geq \ldots 
\geq \lambda_{2^n}$. It follows from the Perron-Frobenius theorem that the 
largest eigenvalue is equal to the spectral norm of $A$, i.e. 
$ \lambda_{\max}(G)=\lambda_1(G)=\Vert A \Vert = \max_{j} |\lambda_j|$.
Note also that for a subgraph of the $n$-cube, or in general, for 
any bipartite graph, $ \lambda_k (G)= - \lambda_{|V|-k}(G),\, k=1,2,
\ldots\,$ and in particular $ |\lambda_{\min}(G)|=\lambda_{\max}(G)$.
It is easy to observe that for every graph $G=(V,E)$ its 
largest eigenvalue $\lambda_1(G)$ is always squeezed between the average
degree of $G$, $\bar{d}=\sum_{v\in V}d_G(v)/|V|$ and its maximal degree
$\Delta(G)=\max_{v\in V}d_G(v)$. In some situtaions these two bounds have the same 
asymptotic value which determines the behavior of the largest eigenvalue. 
For example, this easily gives the asymptotics of the largest egenvalue of
a random subgraph $G(n,p)$ of a complete graph of order $n$ for $p\gg \log n/n$.
On the other hand, in our case there is a 
gap between average and maximal degree of random graph 
$G(Q^n, p)$ for all values of $p<1$ and therefore it is not immediately clear 
how to estimate its largest eigenvalue.

Here we determine the asymptotic value of the largest eigenvalue of
sparse random subgraphs of $n$-cube. To understand better the result, 
observe that if
$\Delta$ denotes a maximal degree of a graph $G$, then $G$ contains a
star $S_\Delta$ and therefore $\lambda_1(G)\ge
\lambda_1(S_\Delta)=\sqrt{\Delta}$. Also, as mentioned above
$\lambda_1(G)$ is at least as large as the average degree of $G$. As for
all values of $p(n)\gg n^{-1}2^{-n}$, a.s. $|E(G(Q^n,p))|=(1+o(1))pn2^n$, we
get that a.s. $\lambda_1(G(Q^n,p))\ge (1+o(1))np$. Combining the above   
lower bounds, we get that a.s. $\lambda_1(G(Q^n,p))\ge
(1+o(1))\max\big(\sqrt{\Delta}, np\big)$. It turns out this lower bound can
be matched by an upper bound of the same asymptotic value, as given by
the following theorem:

\noindent
\begin{theo}
\label{main}
Let $G(Q^n,p)$ be a random subgraph of the $n$-cube and let
$\Delta$ be the maximum degree of $G$. 
Then almost surely the largest eigenvalue of the adjacency matrix
of $G$ satisfies
$$\lambda_1(G)=\big(1+o(1)\big)\max\big(\sqrt{\Delta},np\big),$$
where the $o(1)$ term tends to zero as $\max( \Delta^{1/2}(G), np )$ tends 
to 
infinity.
\end{theo}

As the asymptotic value of the maximal degree of $G(Q^n,p)$ can be easily 
determined for all values of $p(n)$ (see Lemma 2.1), the above theorem
enables us to estimate the asymptotic value of $\lambda_1(G(Q^n,p))$ for all
relevant values of $p$. This theorem is similar to the 
recent result of Krivelevich and Sudakov [11] on the largest eigenvalue of
a random subgraph $G(n,p)$ of a complete graph of order $n$.

The rest of this paper is organized as follows.  
In the next section we gather some necessary technical lemmas about 
random subgraphs of $n$-cube. Section 3 is devoted to the proof of the main 
theorem. Section 4, the last section of the paper, contains some concluding 
remarks.

Throughout the paper we will systematically omit floor and ceiling
signs for the sake of clarity of presentation. All logarithms are 
natural. We will frequently use the inequality $ \left(\begin{array}{c} a \\
b \end{array}\right)
\le
\bigl( ea/b \bigr)^b$. Also we use the following standard notations:
$a_n =\Theta(b_n),\, a_n=O(b_n)$ and $a_n=\Omega(b_n)$
for $a_n >0, b_n >0$ as $n \to \infty$ if there exist
constants $C_1$ and $C_2$ such that $C_1 \* b_n < a_n < C_2 \*
\ b_n,\, a_n < C_2 \* b_n$ or $a_n > C_1 \* b_n$  
respectively. The equivalent notations $a_n=o(b_n)$ and $a_n \ll b_n$
mean that $a_n/b_n \to 0$ as $n \to \infty$. We will say that an event
$\Upsilon_n$, depending on a parameter $n$, holds 
almost surely if 
$\Pr(\Upsilon_n) \to 1 $ as $ n \to \infty $ (please note this phrase, 
very common in the literature on random structures, has a different meaning 
in probability theory, where it means that with probability one 
all but a finite number of events
$ \Upsilon_n$ take place).

\section{Few technical lemmas}
In this section we establish some properties of random subgraphs of the 
$n$-cube which we will use later in the proof of our main theorem.
First we consider the asymptotic behavior of the maximal degree of 
$G(Q^n,p)$. It is not difficult to show that if $p$ is a constant less 
than $1/2$
then a.s. $ \Delta(G)  =(1 +o(1))cn$, where $c$ satisfies the 
equation $ \log 2 + c \log p + (1-c) \log(1-p) = c \log c +(1-c) 
\log(1-c)$ and $ \Delta(G)  =(1 +o(1))n$ for $p\geq 1/2$. 
We omit the proof of this statement since for our purposes it is enough to 
have $ \sqrt{\Delta(G)} =o(np)$ which follows immediately from the fact 
that $\Delta(G) \leq n$. 
The case when $p=o(1)$ is studied in more details in the following lemma.

\begin{lemma}
\label{l2}
Let $G=G(Q^n,p)$ be a random subgraph of the $n$-cube. Denote by
$$\kappa(n)= \max \Big\{ k:  2^n \left(\begin{array}{c} n\\k \end{array}\right)
 p^k (1-p)^{n-k} \geq 1 
\Big\}.$$
Then the following statements hold.

$(i)$\, If $p =o(1)$ and $p$ is not exponentially small in $n$ then
almost surely $\kappa(n)-1 \leq \Delta(G) \leq \kappa(n)+1$.

$(ii)$\, If $p= \Theta( 2^{-n/k} n^{-1})$, then
$2^n \left(\begin{array}{c} n\\k \end{array}\right)
 p^k(1-p)^{n-k} =\Theta(1)$ and  $\kappa(n)= k-1$ or 
$k$. Also almost surely
$\Delta(G)$ is either  $k-1$ or $k$. 

$(iii)$ If $p$ is exponentially small, but not proportional
to $ 2^{-n/k}\* n^{-1},$ then
$\kappa(n) = \bigl [ \frac{n \log 2}{\log (p^{-1}) - \log n} \bigr
]$ and almost surely $\Delta(G)=\kappa(n)$.
\end{lemma}

\noindent
{\bf Proof.}\, Let $X_k$ be the number of vertices of $G(Q^n,p)$ with
degree larger than $k-1$. Then $X_k= \sum_{v \in Q^n} I_v,$
where $I_v$ is  the indicator random variable of the event that
$\deg(v) \geq k$. One can easily calculate the expectation
${\bf E}X_k = 2^n \sum_{l\geq k } 
\left(\begin{array}{c} n\\l \end{array}\right)
p^l (1-p)^{n-l}.$
Also note that  $Q^n$ is bipartite and therefore has independent set
$X$ of size $2^{n-1}$. By definition, the events that $d(v)< k$ are 
mutually independent for all $v \in X$. Therefore we obtain that
\begin{equation}
\label{A}
\Pr\big(X_k=0\big) \leq \prod_{v \in X} \Pr\big(d(v) < k\big)=
\prod_{v \in X}\big(1-{\bf E}I_v\big) \leq \exp\Big(-\sum_{v \in 
X}{\bf E}I_v\Big)=\exp\big(-{\bf E}X_k/2\big).
\end{equation}

Let us now consider the case (i) in more detail. We have that 
$$2^{n-O(np)} \* (n/k)^k p^k  \leq 2^n 
\left(\begin{array}{c} n\\k \end{array}\right)
p^k(1-p)^{n-k} \leq 2^n \* (e \* n / k )^kp^k.$$
Therefore it is easy to check that $\kappa(n)$ must satisfy the 
inequalities
$$\frac{ n\log 2 \big( 1- 1/\log \log (p^{-1})\big)}{\log (p^{-1})}
\leq \kappa(n) \leq
\frac{ n \log 2\big( 1+ 1/\log \log (p^{-1})\big)}{\log (p^{-1})}.$$
By definition, ${\bf E}X_{\kappa(n)}\geq 1$. Elementary computations 
show  that 
${\bf E}X_{k+1} =(1+o(1))\frac{p\log (p^{-1})}{\log 2} {\bf E} X_k$ 
for $k =(1+o(1))\kappa(n)$ which imply  
${\bf E}X_{\kappa(n)-1} \geq (1+o(1)) \frac{\log 2}{p 
\log (p^{-1})}$. Therefore by (1), we have that
\begin{eqnarray*}
\Pr\big(\Delta(G)<\kappa(n)-1\big)&=&\Pr\big(X_{\kappa(n)-1}=0\big) \leq 
\exp\left(-\frac{{\bf E}X_{\kappa(n)-1}}{2}\right)\\
&\leq& 
\exp\left(- (1+o(1))
\frac{\log 2}{p \log(p^{-1})}\right)=o(1).
\end{eqnarray*}
On the other hand, since  ${\bf E}X_{\kappa(n)+1} \leq 1+o(1)$ we have 
that ${\bf E}X_{\kappa(n)+2} \leq (1+o(1)) 
\frac{p\log (p^{-1})}{\log 2}=o(1)$. Thus by Markov's inequality 
we conclude 
that a.s. $X_{\kappa(n)+2}=0$ and therefore 
almost surely $\kappa(n)-1 \leq \Delta(G) \leq \kappa(n)+1$.

Now consider the case (ii). Since 
$p(n)= \Theta( 2^{-n/k}n^{-1})$ we have that
$2^n \left(\begin{array}{c} n\\k \end{array}\right)
p^k(1-p)^{n-k} =\Theta(1)$ and 
$\kappa(n)=k-1$ or $k$. Also it is easy to check that
${\bf E}X_{k+1}=\Theta(2^{-n/k})=o(1)$ and ${\bf E}X_{k-1}=
\Theta(2^{n/k})$. 
Therefore by (1) we have that 
$\Pr\big(\Delta(G)<k-1\big) \leq 
\exp\big(-{\bf E}X_{k-1}/{2}\big)=o(1)$
and by Markov's inequality $\Pr\big(\Delta(G)\geq k+1\big)=o(1)$.

Finally suppose that $p(n)$ is exponentially small but not 
proportional to $2^{-n/k} n^{-1}$ for any $k \geq 1$. Then it is 
rather straightforward to check that 
$\kappa(n) = \bigl [ \frac{n \* \log 2}{\log (p^{-1}) - \log n} 
\bigr ]$ and ${\bf E}X_{\kappa(n)+1} \ll 1 \ll
{\bf E}X_{\kappa(n)}$. Therefore, again using (1) and Markov's 
inequality, we conclude that
$\Pr\big(\Delta(G)<{\kappa(n)}\big)=o(1)$ and 
$\Pr\big(\Delta(G)\geq {\kappa(n)+1}\big)=o(1)$. Thus
almost surely $\Delta(G)=\kappa(n)$. This completes the proof.

Next we need the following lemma, which shows that a.s. $G$ cannot 
have a large number of vertices of high degree too close to one another.
More precisely the following is true. 

\begin{lemma}
\label{l1}
Let $G(Q^n,p)$ be a random subgraph of $n$-cube. Then 
almost surely

\noindent
$(i)$ \, For every  $0 <p \leq 1$ and for any two positive constant $a$ and $b$ such
that $a+b>1$
and $n^b \geq 6np$, $G$ contains no vertex which has within distance
one or two at least $n^a$ vertices of $G$  with degree $ \geq n^b$. 

\noindent
$(ii)$ \, For $p \geq n^{-2/3}$ and any constant $a>0$, $G$ contains no
vertex which has within distance one or two
at least $n^a/p$ vertices of $G$ with degree $ \geq np+np/\log n$.
\end{lemma}

\noindent
{\bf Proof.}\, 
We prove lemma for the case of vertices of distance two, the case 
of vertices of distance one can be treated similarly. 
Note that since the $n$-cube is a bipartite 
graph the vertices which are within the
same distance from a given vertex in $Q^n$ are not adjacent.

${\bf (i)}$ \, Let $X$ be the number of vertices of $G$ which violate condition (i). 
To prove the statement we estimate the expectation of $X$.
Clearly we can choose a vertex $v$ of the $n$-cube in $2^n$ ways. Since there
are at most $n^2$ vertices in $Q^n$ within distance two from $v$ we have at most
$\left(\begin{array}{c} n^2\\n^a \end{array}\right)
$ possibilities to chose a subset $S$ of size $n^a$ of vertices
which will have degree at least $n^b$. 
The probability that the degree of some vertex in $S$ is
at least $n^b$ is bounded by 
$\left(\begin{array}{c} n^2\\n^b \end{array}\right)
p^{n^b}$. Note that these events 
are
mutually independent, since $S$ contains no edges of the $n$-cube. 
Therefore, using that $a+b>1$ and $b>0$, we obtain that
$${\bf E}(X) \leq 2^n \left(\begin{array}{c} n^2\\n^a \end{array}\right)
\left( \left(\begin{array}{c} n^2\\n^b \end{array}\right)
p^{n^b}\right)^{n^a}
 \leq 2^n n^{2n^a}\left(\Big(\frac{enp}{n^b}\Big)^{n^b}\right)^{n^a}
\leq 2^{n+2n^a\log_2n}\, 2^{-n^{a+b}}=o(1).$$ 
Thus by Markov's inequality we conclude that almost surely no vertex 
violates condition (i). 

${\bf (ii)}$ \, Let again $X$ be the number of vertices of $G$ which violate condition
(ii). Similarly as before we have 
$2^n$ choices for vertex $v$ and at most 
$\left(\begin{array}{c} n^2\\n^a/p \end{array}\right)
$ choices for set $S$
of vertices within distance two from $v$ which will have degree $ \geq np+np/\log n$.
Since for all vertices $v\in G$
the degree $d(v)$ is binomially distributed with
parameters $n$ 
and $p$, then by a standard large deviation
inequality  (cf. , e.g., [1], Appendix A)
$${\bf Pr} \Big[d(v) \geq t=np+np/\log n\Big]
\leq e^{-(t-np)^2/2np+(t-np)^3/2(np)^2}=e^{-(1+o(1))np/(2\log^2 n)}.$$
As we already mentioned, the events that vertices in $S$ have 
degree $ \geq np+np/\log n$ are mutually independent.
Therefore, using that $p \geq n^{-2/3}$ and $a>0$, we conclude that
\begin{eqnarray*}
{\bf E}(X) &\leq& 2^n \left(\begin{array}{c} n^2\\n^a /p\end{array}\right)
 \left(e^{-(1+o(1))np/(2\log^2 n)}\right)^
{n^a/p} \leq 2^n n^{2n^a/p} e^{-(1+o(1))n^{1+a}/(2\log^2 n)}\\
&\leq& 
2^n e^{2n^{a+2/3}\log n}e^{-(1+o(1))n^{1+a}/(2\log^2 n)}=o(1).
\end{eqnarray*}
Thus we can complete the proof of the lemma using again Markov's inequality.

\section{Proof of the theorem}
In this section we present our main result. 
We start by listing some simple properties of the largest eigenvalue of
a graph, that we will use later in the proof. Most of these easy statement 
can be found in Chapter 11 of the book of Lov\'asz [12].

\begin{prop}
\label{eigen}
Let $G$ be a graph on $n$ vertices and $m$ edges and with maximum degree
$\Delta$. Let $\lambda_1(G)$ be the largest eigenvalue of the
adjacency matrix of $G$. Then is has the following properties.
\begin{description}
\item[(I)] \label{eigen1}
$\max\big(\sqrt{\Delta}, \frac{2m}{n}\big) \leq
\lambda_1(G) \leq \Delta$.
\item [(II)]
If $E(G)= \cup_i E(G_i)$ then $\lambda_1(G) \leq \sum_i
\lambda_1(G_i)$. If in addition graphs
$G_i$ are vertex disjoint, then $\lambda_1(G)=\max_i \lambda_1(G_i)$
\item[(III)]
If $G$ is a bipartite graph then
$\lambda_1(G) \leq \sqrt{m}$. Moreover if it is a star
of size $\Delta$ then $\lambda_1(G)=\sqrt{\Delta}$.
\item[(IV)]
If the degrees on both sides
of bipartition are
bounded by $\Delta_1$ and $\Delta_2$ respectively, then
$\lambda_1(G) \leq \sqrt{\Delta_1 \Delta_2}$.
\item[(V)] For every vertex $v$ of $G$ let $W_2(G,v)$ denote the number of
walks of length two in $G$ starting at $v$. Then
$\lambda_1(G) \leq \sqrt{\max_v W_2(G,v)}$. 
\end{description}
\end{prop}

\noindent
{\bf Proof of Theorem 1.1.}\,
We already derived in the introduction the lower bound 
of this theorem so we will concentrate on proving an
upper bound. We will frequently use the following simple fact that 
between any two distinct vertices of the $n$-cube there are 
at most two paths of lengths two.
We divide the proof into few cases with respect to the value of the edge
probability $p$. In each case we  partition o $G$ into smaller subgraphs,
whose largest eigenvalue is easier to estimate. 
We start with a rather easy case when the random graph is relatively 
sparse.

{\bf Case 1.}\, Let $e^{-\log^4 n} \leq p \leq n^{-2/3}$. 
For these values of $p$, by Lemma 2.1, we have that 
$\Delta(G) \geq \Omega\big(\frac{n}{\log^4 n}\big)$.
Partition the
vertex set of $G$ into three subsets as follows. Let $V_1$ be the set of
vertices with degree at most $n^{2/5}$, let $V_2$ be the set of vertices
with degree larger than $n^{2/5}$ but smaller than $n^{4/7}$ and let $V_3$
be the set of vertices with degree at least $n^{4/7}$. Also let $G_1$ be a
subgraph of $G$ induced by $V_1$, let $G_2$ be a subgraph induced by
$V_2 \cup V_3$, let $G_3$ be a bipartite graph containing edges of $G$
connecting $V_1$ and $V_2$ and finally let $G_4$ be a bipartite graph
containing edges connecting $V_1$ and $V_3$.
By definition $G=\cup_i G_i$ and thus by claim {\bf (II)} of
Proposition 3.1 we obtain that
$\lambda_1(G) \leq \sum_{i=1}^4 \lambda_1(G_i)$.

Since the maximum degree of graph $G_1$ is at most $n^{2/5}$, then 
by {\bf (I)} it follows that $\lambda_1(G_1) \leq
n^{2/5}$.
The degrees of vertices of bipartite graph $G_3$ are bounded on 
one side by $n^{2/5}$ and on another by $n^{4/7}$. Hence using {\bf (IV)} 
we conclude that 
$\lambda_1(G_3)\leq \sqrt{n^{2/5} n^{4/7}}=n^{17/35}$. 
Let $v$ be a vertex of $G_2$ and let $N_2(G_2,v)$ be the set of vertices of
$G_2$
which are within distance exactly two from $v$. Since between any two
distinct vertices of $Q^n$ there at  
most two paths of length two it is easy to see that the number of walks of
length two
in $G_2$ starting at $v$ is bounded by $d_{G_2}(v)+2|N_2(G_2,v)|$.
Since every vertex of $V_2 \cup V_3$ has degree
in $G$ at least $n^{2/5}$, using Lemma 2.2 (i) with $a=4/5$ and $b=2/5$ we
obtain that almost surely
both $d_{G_2}(v)$ and $|N_2(G_2,v)|$ are bounded by $n^{4/5}$. 
Therefore for every vertex $v$ in $G_2$ we have $W_2(G_2,v)\leq
d_{G_2}(v)+2|N_2(G_2, v)| \leq 3n^{4/5}$ and hence 
by {\bf (V)} $\lambda_1(G_2)\leq \sqrt{3n^{4/5}}=\sqrt{3}n^{2/5}$.

Finally we need to estimate $\lambda_1(G_4)$. To do so consider partition of 
$V_1$ into two parts. Let $V'_1$ be the set of vertices in $V_1$ with at
least two neighbors in $V_3$ and let $V''_1=V_1-V'_1$. 
Let $G'_4$ and $G''_4$ be bipartite graphs with parts $(V'_1,
V_3)$ and $(V''_1,V_3)$ respectively. By definition,
$G_4=G_4' \cup G_4''$ and hence $\lambda_1(G_4) \leq
\lambda_1(G'_4)+\lambda_1(G''_4)$.
Since the vertices in 
$V''_1$ have at most one neighbor in $V_3$ and the graph $G_4''$ is
bipartite it follows that 
$G_4''$ is the union of vertex disjoint stars of size at most
$\Delta(G)$. So
by {\bf (III)} we get $\lambda_1(G''_4)\leq \sqrt{\Delta(G)}$.
Now let $u$ be the vertex of $V_3$ with at least $2n^{1/2}$ neighbors in
$V'_1$. By definition, every neighbor of $u$ in $V'_1$ has also an additional
neighbor in $V_3$, which is distinct from $u$.
Therefore we obtain that there are at least
$2n^{1/2}$ simple paths of length two from $u$ to the set $V_3$. Since
between
any two distinct vertices of the $n$-cube  there are at most two paths of 
length
two we obtain that $u$ has at least $n^{1/2}$ other vertices of $V_3$
within distance two. Since the degree of all these vertices is at least
$n^{4/7}$ it follows from Lemma 2.2 (i) with $a=1/2$ and $b=4/7$ that
a.s. there is no vertex $u$ 
with this property. Therefore the
degree of every vertex from $V_3$ in bipartite graph $G'_4$ is bounded by 
$2n^{1/2}$ and we also have  that the
degree of every vertex from $V'_1$ is at most
$n^{2/5}$. So using again {\bf (IV)} we obtain 
$\lambda_1(G'_4)\leq \sqrt{2n^{1/2}n^{2/5}}=\sqrt{2}n^{9/20}$.
This implies the desired bound on $\lambda_1(G)$, since
$$\lambda_1(G) \leq \sum_i \lambda_1(G_i)\leq
n^{2/5}+\sqrt{3}n^{2/5}+n^{17/35}+
\Big(\sqrt{2}n^{9/20}+\sqrt{\Delta(G)}\Big)=\big(1+o(1)\big)\sqrt{\Delta(G)}.$$

{\bf Case 2.}\, Let $p \geq n^{-4/9}$. This case when the random graph is
dense is also quite simple. Indeed,
partition the vertices of $G$ into two parts. Let $V_1$ be the set of
vertices with degree larger than $np+np/\log n$ and let $V_2$ be the rest of
the vertices. Clearly  $G=\cup_iG_i$, where $G_1$ is a subgraph
induced by $V_1$, $G_2$ is a subgraph induced by $V_2$ and $G_3$ is a
bipartite subgraph
with bipartition $(V_1,V_2)$. By definition, the maximum degree of $G_2$ is
at most $np+np/\log n$, implying $\lambda_1(G_2) \leq np+np/\log n$.
Since  every vertex in $V_1$ has degree at least $np+np/\log n$,
by Lemma 2.2 (ii) with $a=1/18$ we obtain that almost surely no vertex in $G$
can have more than $n^a/p \leq n^{1/2}$ vertices in $V_1$ within distance one
or two. In particular, this implies that the maximum degree in 
$G_1$ is at most $n^{1/2}$ and so $\lambda_1(G_1) \leq n^{1/2}$.

Partition $V_2$ into two parts. Let $V'_2$ be the set of vertices in $V_2$
with at
least two neighbors in $V_1$ and let $V''_2=V_2-V'_2$.
Let $G'_3$ and $G''_3$ be bipartite graphs with parts $(V_1,V'_2)$
and $(V_1,V''_2)$ respectively. By definition,
$G_3=G'_3 \cup G''_3$ and thus $\lambda_1(G_3) \leq
\lambda_1(G'_3)+\lambda_1(G''_3)$.
Since the vertices in
$V''_2$ have at most one neighbor in $V_1$ and the graph $G''_3$ is
bipartite it follows that
$G''_3$ is the union of vertex disjoint stars of size at most
$\Delta(G)$. So
by {\bf (III)} we get $\lambda_1(G''_3)\leq \sqrt{\Delta(G)} \leq n^{1/2}$.
Now let $u$ be the vertex of $V_1$ with at least $2n^{1/2}$ neighbors in
$V'_2$. By definition, every neighbor of $u$ in $V'_2$ has also an additional
neighbor in $V_1$, which is distinct from $u$.
Therefore we obtain that there are at least
$2n^{1/2}$ simple paths of length two from $u$ to the set $V_1$. Since
between
any two distinct vertices there are at most two paths of length
two we obtain that $u$ has at least $n^{1/2}$ other vertices of $V_1$
within distance two. As we already explain in the previous paragraph, this
almost surely does not happen. Thus the
degree of every vertex from $V_1$ in bipartite subgraph $G'_3$ is bounded by
$2n^{1/2}$ and we also have  that the
degree of every vertex from $V'_2$ is at most
$np+np/\log n$. So using  {\bf (IV)} we obtain
$\lambda_1(G'_3)\leq \sqrt{2n^{1/2}(np+np/\log n)}$. Now since $np \geq
n^{5/9}$
it follows that all $\lambda_1(G_1), \lambda_1(G'_3), \lambda_1(G''_3)=o(np)$
and therefore 
$$\lambda_1(G) \leq \sum_i \lambda_1(G_i)\leq \lambda_1(G_2)+o(np) \leq 
np+np/\log n+o(np)=(1+o(1))np.$$

{\bf Case 3.}\, Let $n^{-2/3} \leq p \leq n^{-4/9}$. This part of the proof  
is slightly more involved then two previous ones since in particular it 
needs to deal with a delicate case when $np$ and $\sqrt{\Delta(G)}$ are
nearly equal.

Partition the vertex set of $G$ into four parts.
Let $V_1$ the set of vertices with degree at least
$n^{2/3}$ and let $V_2$ be the set of vertices with degrees larger
than $np+np/\log n$ but less than $n^{2/3}$.
Let $V_4$ contains all vertices $G$ which have at least one neighbor in $V_1$
and degree at most $np+np/\log n$.
Finally let $V_3$ be the set of remaining vertices of $G$.
Note that by definition there are no edges between
$V_1$ and $V_3$ and every vertex from $V_3$ also have degree at most $np+np/\log n$ 
in $G$.

We consider the following subgraphs of $G$.
Let $G_1$ be the bipartite subgraph 
containing all the edges between $V_1$ and $V_4$. 
Partition $V_4$ into two parts. Let $V'_4$ be the set of vertices in $V_4$  
with at
least two neighbors in $V_1$ and let $V''_4=V_4-V'_4$.
Let $G'_1$ and $G''_1$ be bipartite graphs with parts $(V_1,V'_4)$
and $(V_1,V''_4)$ respectively. By definition,
$G_1=G'_1 \cup G''_1$ and thus $\lambda_1(G_1) \leq
\lambda_1(G'_1)+\lambda_1(G''_1)$.
Since the vertices in
$V''_4$ have at most one neighbor in $V_1$ and the graph $G''_1$ is
bipartite it follows that
$G''_1$ is the union of vertex disjoint stars of size at most
$\Delta(G)$. So
by {\bf (III)} we get $\lambda_1(G''_1)\leq \sqrt{\Delta(G)}$.
Now let $u$ be the vertex of $V_1$ with at least $2n^{2/5}$ neighbors in
$V'_4$. By definition, every neighbor of $u$ in $V'_4$ has also an additional
neighbor in $V_1$, which is distinct from $u$.
Therefore we obtain that there are at least 
$2n^{2/5}$ simple paths of length two from $u$ to the set $V_1$. 
Similar as before, this implies that 
$u$ has at least $n^{2/5}$ other vertices of $V_1$
within distance two. By Lemma 2.2 (i) with $a=2/5$ and $b=2/3$ this
almost surely does not happen. Thus the
degree of every vertex from $V_1$ in bipartite subgraph $G'_1$ is bounded by
$2n^{2/5}$ and we also have  that the
degree of every vertex from $V'_4$ is at most
$np+np/\log n \leq n^{5/9}$. So using  {\bf (IV)} we obtain
$\lambda_1(G'_1)\leq
\sqrt{2n^{2/5}n^{5/9}}=\sqrt{2}n^{43/90}=o(\sqrt{\Delta})$. Therefore 
$\lambda_1(G_1) \leq \lambda_1(G'_1)+\lambda_1(G''_1) \leq 
(1+o(1))\sqrt{\Delta}$.

Our second subgraph $G_2$ is induced by the set $V_3$.
By definition, the maximum degree in it is at most 
$np+np/\log n$ and therefore $\lambda_1(G_2) \leq (1+o(1))np$.
Crucially this graph is vertex disjoint form $G_1$ which implies by 
{\bf (II)} that 
$$\lambda_1\big(G_1 \cup
G_2\big)=\max\big(\lambda_1(G_1),\lambda_1(G_2)\big)\leq 
\big(1+o(1)\big)\max\left(\sqrt{\Delta},np\right).$$

Next we define the remaining graphs whose union with $G_1$ and $G_2$ equals
to $G$ and show that
their largest eigenvalues contribute only smaller order terms in the 
upper bound on $\lambda_1(G)$. Let $G_3$ be the subgraph of $G$ induced by
the set $V_1 \cup V_2$. By definition, every vertex in $G_3$ have at least 
$np+np/\log n$ neighbors in $G$. Therefore by Lemma 2.2 (ii) with $a=1/12$ we
obtain that for every $v \in G_3$ there at most $n^a/p \leq n^{3/4}$ other 
vertices of $G_3$ within distance one or two. This implies that
$d_{G_3}(v)$ and $|N_2(G_3,v)|$ are both bounded by $n^{3/4}$. Then, as we already 
show in Case 1, the total number of walks of length two starting at $v$ is
bounded by $d_{G_3}(v)+2|N_2(G_3,v)| \leq 3n^{3/4}$. Thus 
by {\bf (V)} we get $\lambda_1(G_3) \leq \sqrt{3n^{3/4}}=o(\sqrt{\Delta})$.
  
Let $u$ be a vertex of $V_3 \cup V_4$ which has at least $2n^{2/5}$
neighbors in the set $V_4$. Since every vertex in $V_4$ have at least
one neighbor in $V_1$ we obtain that there at least $2n^{2/5}$
simple paths of length two from $u$ to $V_1$. On the other hand we know that
there are at most two such paths between any pair of distinct vertices.
This implies that $u$ has at least $n^{2/5}$ vertices within distance two 
whose degree is at least $n^{2/3}$. Using Lemma 2.2 (i) with $a=2/5$ and
$b=2/3$ we conclude that almost surely there is no such vertex $u$.
Now let $G_4$ be a subgraph induced by the set $V_4$ and let
$G_5$ be the bipartite graph with parts $(V_3,V_4)$. By the above discussion,
the maximum degree of $G_4$ is at most $2n^{2/5}$, implying
$\lambda_1(G_4) \leq 2n^{2/5} =o(\sqrt{\Delta})$. 
We also know that every vertex from
$V_3$ has at most $2n^{2/5}$ neighbors in $V_4$ and every vertex in $V_4$
have at most $np+np/\log n \leq n^{5/9}$ neighbors in $V_3$. Therefore by
{\bf (IV)} we obtain that $\lambda_1(G_5) \leq 
\sqrt{2n^{2/5}n^{5/9}}=\sqrt{2}n^{43/90}=o(\sqrt{\Delta})$.

Finally consider the bipartite subgraph $G_6$ whose parts are $V_2$ and 
$V_3 \cup V_4$. Let $X$ be the set of vertices from $V_3 \cup V_4$ with
at least $2n^{2/7}$ neighbors in $V_2$ and let $Y=V_3 \cup V_4-X$.
Note that $G_6=G'_6 \cup G''_6$ where $G'_6$ is bipartite graph with parts 
$(V_2,X)$ and $G''_6$ is bipartite graph with parts $(V_2,Y)$.
The upper bound on $\lambda_1(G''_6)$ follows immediately from the facts 
that $G''_6$ is bipartite,
the degree of vertices in $V_2$ is bounded by $n^{2/3}$ and, by definition, 
every vertex in $Y$ has at most $2n^{2/7}$ neighbors in $V_2$.
Therefore $\lambda_1(G''_6) \leq \sqrt{2n^{2/7}n^{2/3}}=
\sqrt{2}n^{10/21}=o(\sqrt{\Delta})$. To bound 
$\lambda_1(G'_6)$, note that almost surely every vertex in $V_2$
has at most $n^{3/7}$ neighbors in $X$. Indeed, let $u \in V_2$ be the
vertex
with more than $n^{3/7}$ neighbors in $X$. 
Since every neighbor of $u$ in $X$ has at least
$2n^{2/7}-1$ additional neighbors in $V_2$ different form $u$ we obtain that there 
at
least $(2n^{2/7}-1)n^{3/7}=(2+o(1))n^{5/7}$
simple paths of length two from $u$ to $V_2$. On the other hand we know that
there are at most two such paths between any pair of distinct vertices.
This implies that $u$ has at least $(1+o(1))n^{5/7}$ vertices of $V_2$ within
distance two. By definition, every vertex of $V_2$ has at
least $np+np/\log n$ neighbors in $G$. Therefore 
using Lemma 2.2 (ii) with $a=1/22$
we conclude that almost surely there is no such vertex $u$.
Now the upper bound on $\lambda_1(G'_6)$ can be obtained using that
$G'_6$ is bipartite,
the degree of vertices in $X$ is bounded by $np+np/\log n \leq n^{5/9}$ and
that every vertex in $V_2$ has at most $2n^{3/7}$ neighbors in $X$.
Indeed, by {\bf (IV)} $\lambda_1(G'_6) \leq \sqrt{2n^{3/7}n^{5/9}}=
\sqrt{2}n^{31/63}=o(\sqrt{\Delta})$ and hence
$\lambda_1(G_6) \leq \lambda_1(G'_6)+\lambda_1(G''_6)=o(\sqrt{\Delta})$.

From the above definitions it is easy to check that $G=\cup_{i=1}^6 G_i$.
Hence using our estimates on the largest eigenvalues of graphs $G_i$ we 
obtain the desired upper bound on $\lambda_1(G)$, as follows
\begin{eqnarray*}
\lambda_1(G) &\leq& \lambda_1\big(G_1 \cup G_2\big)+\sum_{i \geq 3}
\lambda_1(G_i) \leq 
(1+o(1))\max\left(\sqrt{\Delta(G)}, np\right)+o\big(\sqrt{\Delta(G)}\big)\\
&=&
(1+o(1))\max\left(\sqrt{\Delta(G)}, np\right).
\end{eqnarray*}
This completes the proof of the third case. Now to finish the proof of the 
theorem it remains to deal with the last very simple case when the 
random graph is very sparse.

{\bf Case 4.}\, Let $p \leq e^{-\log^4n} $.  For every integer $k \geq 1$ 
denote by $Y_k$ the number of connected components with 
$k$ edges. It is not difficult to see that
${\bf E}Y_k \leq O(2^n k! n^k p^k)$. Indeed, we can pick the first vertex 
in the connected component in $2^n$ ways. Suppose we already know the first 
$1 \leq t \leq k$ vertices of the component. Then these vertices are 
incident to at most $tn$ edges of the $n$-cube and therefore we can pick the 
next edge only in at most $tn$ ways. This gives at most $2^n \prod_{t=1}^k tn=
2^n k! n^k$ ways to pick the edges of the connected component.

First consider the case when $p$ is not exponentially small. 
Then, by Lemma 2.1 we have that almost surely
$\Delta(G)=(1+o(1))\kappa(n)$, where 
$\kappa(n)= \max \{ k:  2^n \left(\begin{array}{c} n\\k \end{array}\right)
 p^k (1-p)^{n-k} \geq 1 \}$ and 
$\kappa(n)$ tends to infinity together with $n$.
Let $k_0=\kappa(n)+\kappa(n)/\log \kappa(n)$. Then it is easy to 
check that ${\bf E}Y_{k_0}=o(1)$ and therefore, by Markov's inequality, 
almost surely $G(Q^n,p)$ contains no connected component with more than $k_0$ 
edges. Since the largest eigenvalue of $G$ is the maximum of the eigenvalues 
of its connected components and  the largest eigenvalue of a component with 
$k_0$ edges is not greater than $\sqrt{k_0}$  (see, parts {\bf (II)} and 
{\bf (III)} of Proposition 3.1), we obtain that
$$\lambda_1(G)\leq \sqrt{k_0} \leq \sqrt{\kappa(n)}+
\sqrt{\kappa(n)/\log \kappa(n)}=(1+o(1))\sqrt{\Delta(G)}.$$

Next, let  $p \leq 2^{-\alpha n}$ for some fixed $\alpha>0$.
If $p$ is not proportional to  $2^{-n/k}n^{-1},\, k=1,2,3, \ldots,$ 
then it follows from part (iii) of  Lemma 2.1 that with
probability going to one the maximum degree of $G(Q^n,p)$ is 
$\kappa(n)=\bigl [ \frac{n \* \log 2}{\log (p^{-1}) - \log n} \bigr ]$. 
Note that in this case $\kappa(n)$ is a constant and it is easy to check 
that ${\bf E}Y_{\kappa(n)+1} \leq 
O\big(2^nn^{\kappa(n)+1}p^{\kappa(n)+1}\big)=o(1)$.
Thus, by Markov's inequality, there are no connected components with more 
than $\kappa(n)$ edges.
Since the largest eigenvalue of $G$ is the maximum of the eigenvalues of 
its connected components and  the largest eigenvalue of a component with $k$ 
edges is not greater than $\sqrt{k}$ ( and is equal to 
$\sqrt{k}$ only if the component is a star on $k+1$ vertices),
we obtain that a.s. $\lambda_1(G)=\sqrt{\kappa(n)}=\sqrt{\Delta(G)}.$

  Finally if  $p(n)$ is proportional to  $2^{-n/k} n^{-1}, \, k=1,2,3,
\ldots, $ then by part (ii) of Lemma 2.1 almost surely $\Delta(G)  \in 
\{ 
k-1,k\} \ $ and again one can check that ${\bf E}Y_{k+1}$ is exponentially 
small. Using Markov's inequality, as before, we conclude that there are no
connected components with more than $k$ edges.
Therefore a.s. $\lambda_1(G)$ is either $\sqrt{\Delta(G)}$ or 
$\sqrt{\Delta(G)+1}$. This completes the proof of the theorem.

\section{Concluding remarks}

 There are several other important questions that are beyond the reach of 
the presented technique. The most fundamental is perhaps the local 
statistics
of the eigenvalues, in particular the local statistics near the edge of the 
spectrum. For results in this direction for other random matrix models
we refer the reader  to [17], [18], [16].
A recent result of Alon, Krivelevich and Vu [3]
states that the deviation of the first, second, etc. largest 
eigenvalues from its mean is at most of order of $O(1)$. Unfortunately
our results  give only the leading term of the mean.

 A second, perhaps even more difficult question  is whether the local 
behavior of the eigenvalues is sensitive to the details of the distribution
of the matrix entries of $ \ A. \ $
We refer the reader to
[15], [6], [16], [10] for
the results of that nature for unitary invariant and Wigner random matrices.

\vspace{0.25cm}
\noindent
{\bf Acknowledgements.}\,  The first author would like to thank Sergey 
Gavrilets and Janko Gravner for bringing this problem to his attention 
and for useful discussions.

\def\cmp{{\it Commun. Math. Phys.} }

\end{document}